\input amstex
\input Amstex-document.sty
\widestnumber\key{41}
\nologo
\TagsAsMath
\redefine\div{\operatorname{div}}
\define\osc#1{\operatorname{osc} (#1)}
\define\Ah{$\Cal A$-harmonic}
\define\sh{$\Cal A$-superharmonic}

\define\A{\Cal A}
\define\rn{\bold R^n}
\define\re{\bold R}
\redefine\O{\Omega}
\define\dome{\partial\O}

\define\xo{x_0}
\define\kapu{\operatorname{cap}_{p}}
\define\vep{\varepsilon}
\define\loc{\operatorname{loc}}
\define\wo#1{W^{1,p}(#1)}
\define\lwo#1{W^{1,p}_{\loc}(#1)}
\define\woo#1{W_0^{1,p}(#1)}

\define\lwome{W^{1,p}_{\loc}(\O)}
\define\woome{\woo{\O}}
\define\caa#1{C^{\infty}(#1)}
\define\cao#1{C_0^{\infty}(#1)}
\define\vf{\varphi}
\define\wL#1#2{\operatorname{weak}-L^{#1}(#2)}
\define\wLl#1#2{\operatorname{weak}-L^{#1}_{\loc}(#2)}
\font\boldmathHuge = cmmib10 scaled 1748
\font\boldmathlarge = cmmib10 scaled 1200

\pageno 167

\def\shorttitle{$p$-Laplacian Type Equations Involving Measures}
\def\shortauthor{T. Kilpel\"ainen}

\topmatter %
\title\nofrills{\boldHuge {\boldmathHuge p}-Laplacian Type Equations Involving Measures}
\endtitle

\author \Large T. Kilpel\"ainen* \endauthor

\thanks *Department of Mathematics \& Statistics, University of
Jyv\"askyl\"a, PO Box 35, FIN-40014 Jyv\"askyl\"a, Finland.
E-mail: terok\@math.jyu.fi \endthanks

\abstract\nofrills \centerline{\boldnormal Abstract}

\vskip 4.5mm

{\ninepoint This is a survey on problems involving equations $-\operatorname{div}{\Cal A}(x,\nabla u)=\mu$, where
$\mu$ is a Radon measure and ${\Cal A}:\bold {R}^n\times\bold R^n\to \bold R^n$ verifies Leray-Lions type
conditions. We shall discuss a potential theoretic approach when the measure is nonnegative. Existence and
uniqueness, and different concepts of solutions are discussed for general signed measures.

\vskip 4.5mm


}
\endabstract
\endtopmatter

\document

\baselineskip 4.5mm \parindent 8mm

\specialhead \noindent \boldLARGE 1. Introduction \endspecialhead

Throughout this paper we
let $\O$ be an open set in $\rn$ and  $1<p<\infty$ a fixed number.
We shall consider equations
$$
 -\div\A(x,\nabla u)=\mu\,,\tag{1.1}
$$
where $\mu$ is a Radon measure. We suppose that the mapping
$\A\:\rn\times\rn\to\rn$, $(x,\xi)\mapsto \A(x,\xi)$,
is measurable in $x$ and continuous in $\xi$ and that it verifies
the structural conditions:
$$
\gathered
\Cal A(x,\xi)\cdot \xi\,\ge\lambda|\xi|^p  \,,  
\quad
|\Cal A(x,\xi)|\,\le\,\Lambda|\xi|^{p-1}, \quad\text{ and }\\
\left(\Cal A(x,\xi)-\Cal A(x,\zeta)\right)\cdot(\xi-\zeta)\,>\,0\,,
\endgathered
\tag{1.2}
$$
for  a.e\. $x\in \rn$ and all $\xi\ne \zeta\in\rn$.
A prime example of the operators is the $p$-Laplacian
$$
- \Delta_p u=-\div(|\nabla u|^{p-2}\nabla u).
$$

In Section 2 we discuss how nonlinear potential theory is
related to equations like (1.1); it corresponds to nonnegative
measures.
Then in
Section 3
we discuss the existence and uniqueness for  (1.1)
with general measures.

\specialhead \noindent \boldLARGE 2. Potential theoretic approach
\endspecialhead

A continuous solution $u\in\lwome$ of
$ -\div\A(x,\nabla u)=0$
is called $\A$-{\it harmonic} in $\O$. An  $\A$-{\it superharmonic}
function  in $\O$ is a lower semicontinuous function
$u\:\O\to\re\cup\{\infty\}$ that is not identically $\infty$ in any
component of $\O$ and that obeys the following comparison property:
{\sl for each open $D\subset\subset\O$ and each $h\in C(\bar D)$,
\Ah\ in $D$, the inequality $u\ge h$ on $\partial D$ implies
$u\ge h$ in $D$.}

For $k>0$ and $s\in\re$, let
$$
  T_k(s)=\max\big(-k,\min(s,k)\big)
$$
be the truncation operator. Then we have

{\bf 2.1. Theorem.} \cite{12, 26, 15}
{\it If $u$ is \sh\ in $\O$, then
$$
 T_k(u)\in\lwome\quad\text{ for all }k>0\,.
$$}

This enables us to show  that \sh\ functions have a ``gradient'':
Suppose that a function $u$ that is finite a.e\. has the property that
$T_k(u)\in\lwome$ for all $k>0$. Then we define the (weak) {\it gradient}
of $u$ as
$$
\nabla u(x)= DT_k(u)\quad\text{ if } |u(x)|<k\,.
$$
Here $ DT_k(u)$ is the distributional gradient of the Sobolev
function $T_k(u)\in\lwome$. Then $\nabla u$ is well defined. Observe that
if $\nabla u$ is locally integrable, then it is the distributional
derivative of $u$. However, it may happen that $u$ or $\nabla u$ fails to
be locally integrable and so $\nabla u$ is not always distributional
derivative, see \cite{15}, \cite{8}; this is a real issue only for
$p<2-1/n$.

{\bf 2.2. Theorem.} \cite{26, 12, 15, 10, 1}
{\it Suppose that $u$ is \sh\ in $\O$.
\roster
\item"i)" If $1<p<n$, then
$$
u\in \wLl{n(p-1)/(n-p)}{\O}\quad \text{ and }\quad \nabla u\in
\wLl{n(p-1)/(n-1)}{\O}.
$$
\item"ii)" If $p=n$, then $u$ is locally in $BMO$ and hence
$u\in L^q_{\operatorname{loc}}(\O)$ for all $q>0$;
moreover
$$
\nabla u\in \wLl{n}{\O}.
$$
 \item"iii)" If $p>n$, then
$$
u\in \lwome
$$
and hence $u$ is (H\"older) continuous.
\endroster
}

If $p>n$ \sh\ functions are continuous and locally in $W^{1,p}$, moreover
Radon measures then are in the dual of Sobolev space $W^{1,p}$. This makes
the cases $p>n$ very special and quite simple by the classical
results of Leray and Lions. Henceforth we shall be
concerned mainly with cases $1<p\le n$.

For \sh\ $u$ we have by Theorem
2.2 that $|\nabla u|^{p-1}$ is locally integrable. So the
distribution
$$
-\div\A(x,\nabla u)(\vf):=\int_{\O}\A(x,\nabla u)\cdot\nabla\vf\,dx
\,,\quad\vf\in \cao{\O}\,,
$$
is well defined.
\sh\ functions give rise to equation (1.1).

{\bf 2.3. Theorem.} \cite{18}
{\it If $u$ is \sh\ in $\O$, then $-\div\A(x,\nabla u) $
is represented by a nonnegative Radon
measure $\mu$.
}

As to the existence we have:

{\bf 2.4. Theorem.} \cite{18, 3}
{\it Given a finite nonnegative Radon measure $\mu$ on bounded $\O$,
there is an \sh\
function $u$ in $\O$ such that
$$
\cases
-\div\A(x,\nabla u)=\mu\\
T_k(u)\in\woome\text{ for all }k>0\,.
\endcases\tag{2.5}
$$
}

If $\O$ is undounded, then there also is an \sh\ solution to (1.1).
This is rather easily seen if $1<p<n$. The case $p\ge n$ requires a
more careful analysis, see \cite{9, 17}.

In light of Theorem 2.2 we have that a solution $u$ to
(2.5) satisfies
$$
u\in W^{1,q}_0(\O)\quad\text{ for all } q<\frac{n(p-1)}{n-1}
$$
(for $p>2-1/n$).
One naturally asks if such a solution is unique. Unfortunately this is
not the case in general (except for $p>n$).
To see this consider the function
$$
 u(x)=\cases |x|^{\frac{p-n}{p-1}} -1&\text{ if }1<p<n ,\\
  -\log|x|&\text{ if }p=n\,.
\endcases
$$
In the $p$-Laplacian cases,  $u$ is then
a $p$-superharmonic solution to (2.5)
with $\mu$=0 in the punctured
ball $\O=B(0,1)\setminus\{0\}$, but $v\equiv 0$ is
another solution. This is rather artificial example but there are more
severe ones. The question
of uniqueness is a real issue to which we shall return in Section
3 below.

\specialhead \noindent \boldlarge Regularity and estimates
\endspecialhead

In the classical potential theory
the uniqueness can be solved by the aid of the Riesz decomposition theorem
which states that superharmonic functions are sums of a potential and
a harmonic function. No such decomposition is available in the nonlinear
world. However, this lack can be compensated for to an extent by estimating
in terms of the {\it Wolff potential} of the measure $\mu$,
$$
{\bold W}^{1,p}_{\mu}(x_0,r)=
\int_0^r\left(\frac{\mu(B(x_0,t))}{t^{n-p}}\right)^{\frac 1{p-1}}\frac{dt}{t}\,.
$$

{\bf 2.6. Theorem.} \cite{19}
{\it Let $u$ be a nonnegative \sh\ function in $B(\xo,3r)$ with
$\mu=-\div\A(x,\nabla u)$.
Then
$$
 c_1{\bold W}^{1,p}_{\mu}(x_0,r)
\le u(x_0)\le  c_2{\bold W}^{1,p}_{\mu}(x_0,2r)
 + c_3\inf_{B(x_0,r)}u\,,
$$
where $c_j=c_j(n,p,\lambda,\Lambda)>0$.
}

Theorem 2.6 was discovered by the author with Mal\'y
\cite{18, 19} and later generalized for equations depending also on
$u$  by Mal\'y and Ziemer \cite{30}. Mikkonen \cite{32} worked out
the argument for weighted operators, this was later written up in a
metric space setup in \cite{2}.
Recently a totally different proof that works for quasilinear subelliptic
operators was found by Trudinger and Wang \cite{39}. Labutin \cite{22}
gave a generalization for $k$-Hessian operators.

As the first major application
of the potential estimate in 2.6 the author
and Mal\'y
established the necessity of the Wiener test for the regularity
for the Dirichlet problem: We say that $\xo\in\dome$
is an $\A$-{\it regular} boundary point of bounded $\O$ if
$$
\lim_{x\to\xo}u(x)=\vf(\xo)
$$
whenever $\vf\in\caa{\rn}$ and $u$ is \Ah\ in $\O$ with
$u-\vf\in\woome$. Then it turns out that regularity is independent
of the particular operator and it depends only on its type $p$.
More precisely,
define the $p$-{\it capacity} of the set $E$ as
$$
\kapu(E):=\inf\int_{\rn}|v|^p+|\nabla v|^p\,dx\,,
$$
where the infimum is taken over all
$v\in W^{1,p}(\rn)$ such that
$v\ge 1$  on an open neighborhood of  $E$. Then

{\bf 2.7. Theorem.} \cite{31, 19}
{\it A boundary point $\xo\in\dome$ is regular if and only if
$$
\int_0^1\left(\frac{\kapu(\complement\O\cap
B(x_0,t))}{t^{n-p}}\right)^{\frac 1{p-1}}
\frac{dt}{t}=\infty\,.
$$
}

Maz'ya \cite{31} introduced the Wiener type test in 2.7 and
proved its sufficiency. Gariepy and Ziemer \cite{11} generalized the result
by a different argument. Lindqvist and Martio \cite{27}
proved the necessity for $p>n-1$; the general case was treated in \cite{19}.
For generalizations see \cite{30, 32, 39, 22}.

Theorem 2.6 can also be used to characterize singular solutions
and (H\"older) continuity of \sh\ functions (see
\cite{17, 19, 21}). Recall that \Ah\ functions are locally
H\"older continuous: there is a constant $\varkappa\in (0,1]$ depending
only on the structure such that
$$
\osc{u,B(x,r)}\le c(\frac rR)^{\varkappa}
\osc{u,B(x,R)}\tag{2.8}
$$
whenever $u$ is \Ah\ in $B(x,R)$, $r<R$ \cite{37, 15}.
For the $p$-Laplacian $\varkappa=1$.
We have

 {\bf 2.9. Theorem.} \cite{21}
{\it Suppose that $u\in\lwome$ is \sh\ with $\mu=-\div\A(x,\nabla u)$. If $u$
is $\alpha$-H\"older continuous, then
$$
 \mu\big(B(x,r)\big)\le c\,r^{n-p+\alpha(p-1)}\quad
\text{ whenever }\quad B(x,2r)\subset\O.
$$
The converse holds if $\alpha<\varkappa$.
}

The case $\alpha=\varkappa$
is quite different and it is not yet well understood.

Rakotoson and Ziemer \cite{36} proved that the condition of Theorem
2.9 for the measure gives the H\"older continuity of
the solution with some exponent. Lieberman \cite{25} showed that
for smooth operators like the $p$-Laplacian  the condition
$ \mu\big(B(x,r)\big)\le c\,r^{n-1+\vep}$ for some $\vep>0$
implies that the solution is in $C^{1,\beta}$.

Theorem 2.9 can be employed to establish the following
removability result which is due to Carleson \cite{5} in the Laplacian
case.

{\bf 2.10. Theorem.} \cite{21}
{\it Suppose that $u\in C^{0,\alpha}(\O)$ is \Ah\ in $\O\setminus E$.
If $E$ is of $n-p+\alpha(p-1)$ Hausdorff measure zero, then $u$ is \Ah\
in $\O$.

If $E$ is of positive $n-p+\alpha(p-1)$ Hausdorff measure and
$\alpha<\varkappa$, then there is $u\in C^{0,\alpha}(\O)$ that is \Ah\ in
$\O\setminus E$ but not
in the whole $\O$.
}

\specialhead \noindent \boldLARGE 3. General Radon measures
\endspecialhead

In this section we let $\mu$ be any signed Radon measure and consider equation
(1.1). More specifically, we shall discuss the problem
$$
\cases
-\div\A(x,\nabla u)=\mu\\
u=0\text{ on } \dome\,
\endcases\tag{3.1}
$$
where $\O$ is a bounded domain in $\rn$. Here the equation is understood in the
distributional sense, i.e\.
$$
 \int_{\O}\A(x,\nabla u)\cdot\nabla\vf\,dx=\int_{\O}\vf\,d\mu\,,\quad\vf\in\cao{\O},
$$
where we, of course,
assume that $x\mapsto \A(x,\nabla u)$ is locally integrable.
The boundary values $u=0$ are assumed in a weak Sobolev space sense.
The existence of the solution to this problem is known:

{\bf 3.2. Theorem.} \cite{3, 8, 9, 10}
{\it For each Radon measure $\mu$ of finite total variation,
there is a solution $u$ to (3.1) such that
\roster
\item"i)" the truncations
$$
T_k(u)\in\woome\text{ for all }k>0\,,
$$
\item"ii)" $$
u\in \wL{n(p-1)/(n-p)}{\O}\text{ if }1<p<n \quad
\text{ and }\quad
u\in\operatorname{BMO}\text{ if }p=n\,,
$$
\item"iii)"
 $$
\nabla u\in \wL{n(p-1)/(n-1)}{\O}.
$$
\endroster
}

We next discuss the uniqueness of such a solution.
There are examples of mappings $\A$ for which there is a solution $u$
of equation
$\A(x,\nabla u)=0$
such that $u$ satisfies ii) and iii) of Theorem 3.2, but
fails to be \Ah, see \cite{33, 38, 16, 29}.
See also the nonuniqueness example in an irregular domain
after Theorem 2.4 above.

There are various approaches trying to treat the uniqueness problem
by attaching additional attributes to the solution. To formulate these
we need to recall a decomposition of measures.
The $p$-capacity, defined in Section 2, is an outer measure.
Hence the usual proof of the Lebesgue decomposition theorem gives us
that any Radon measure $\mu$ can be decomposed as
$$
  \mu =\mu_0+\mu_s\,,
$$
where $\mu_0$ and $\mu_s$ are Radon measures such that
$\mu_0$ is {\it absolutely continuous} with respect to the
$p$-capacity (i.e., $\mu_0(E)=0$ whenever $\kapu(E)=0$) and
$\mu_s$ is {\it singular} with respect to the
$p$-capacity (i.e., there is a Borel set $B$ such that  $\kapu(B)=0$
and $\mu_s(E\setminus B)=0$ for all $E$).

Let  $u$ be a solution to (3.1) described in Theorem
3.2. We say that
\roster
\item"-" $u$ is an {\it entropy solution} of (3.1) if $u$ is Borel
measurable and
$$
\int_{\O}\A(x,\nabla u)\cdot\nabla T_k(u-\vf)\,dx
\le\int_{\O}  T_k(u-\vf) \,d\mu
$$
for all $\vf\in\cao{\O}$ and $k>0$.
\item"-"  $u$ is a {\it renormalized solution} of (3.1)
if for all $h\in W^{1,\infty}(\re)$ such that $h'$ has
compact support we have
$$
\align
\int_{\O}&\A(x,\nabla u)\cdot\nabla u h'(u)\vf\,dx
+\int_{\O}\A(x,\nabla u)\cdot h(u)\nabla \vf\,dx
\\&=\int_{\O}  h(u)\vf \,d\mu_0
+h(\infty)\int_{\O}  \vf \,d\mu_s^+-
h(-\infty)\int_{\O}  \vf \,d\mu_s^-
\endalign
$$
whenever
$\vf\in W^{1,r}(\O)\cap L^{\infty}(\O)$ with $r>n$ is such that
  $h(u)\vf\in \woome$;
here
$$
h(\infty)=\lim_{t\to\infty} h(t)\,,\quad
h(-\infty)=\lim_{t\to-\infty} h(t)\,,
$$
and
$\mu_s^+$ and $\mu_s^-$ are the positive and negative parts
of the singular measure $\mu_s$.
\endroster

Observe that we assume here that $u$ satisfies equation (3.1)
in the distributional sense. The entropy condition in this context was
first used in \cite{1}; the renormalized solution was introduced by
Lions and Murat \cite{28} and in the refined form in \cite{7, 8}.
The artificial function we had as a counterexample for the uniqueness
(after Theorem 2.4) is not entropy nor renormalized solution.
Existence is known; most existence proofs follow a similar idea to that
in \cite{3}. The uniqueness can be established for certain measures since
one can use truncated solutions as test functions.

{\bf 3.3. Theorem.} \cite{4, 20, 8, 40}
{\it Suppose that $\mu$ is a finite Radon measure. Then
there is a renormalized and an entropy solution of
(3.1). Moreover such a solution is unique if
$\mu$ is absolutely continuous with respect to the $p$-capacity.
}

The uniqueness with $\mu\in L^1$ was proved in \cite{1} and \cite{28}
see also \cite{33}, \cite{34}. For measures absolutely continuous
wrt $p$-capacity, the uniqueness is established e.g. in
 \cite{4, 20, 8, 40}.

A renormalized solution is always an entropy solution,
whence the concepts coincide at least if $\mu$ is
absolutely continuous wrt $p$-capacity; see \cite{8, 6}.

In case of a general measure
the uniqueness of renormalized solution appears to be an open
problem. There are some partial results: Assume that the following strong
monotoneity assumption holds. For all $0\not=\xi,\eta\in\rn$
$$
\big(\A(x,\xi)-\A(x,\eta)\big)\cdot(\xi-\eta)
\ge\cases
\beta|\xi-\eta|^p&\text{ if }p\ge 2,\\
\beta\dfrac{|\xi-\eta|^2}{(|\xi|+|\eta|)^{2-p}}&\text{ if }p< 2.
\endcases\tag{3.4}
$$
Assume also the H\"older continuity:
$$
|\A(x,\xi)-\A(x,\eta)|
\le\cases
\gamma\big(b(x)+|\xi|+|\eta|\big)^{p-2}|\xi-\eta|^2&\text{ if }p\ge 2,\\
\gamma |\xi-\eta|^{p-1} &\text{ if }p< 2\,,
\endcases\tag{3.5}
$$
where $b\in L^p$ is nonnegative. For instance,
the $p$-Laplacian satisfies these assumptions.
Then we have:

{\bf 3.6. Theorem.} \cite{7, 8, 14}
{\it Suppose the additional assumptions (3.4) and
(3.5) hold. If $u$ and $v$ are two renormalized solutions of
(3.1) with measure $\mu$ such that either
$\nabla u-\nabla v\in L^p(\O)$ or $u-v$ is bounded from one side,
then $u=v$.
}

Rakotoson \cite{35} proved that a continuous renormalized solution is
unique in smooth domains; continuity requires the measure be rather special.

\specialhead \noindent \boldlarge Borderline case {\boldmathlarge p} = {\boldmathlarge n}
\endspecialhead

We close this paper by considering the special case when $p=n$.
Then the uniqueness can be reached:

{\bf 3.7. Theorem.} \cite{41, 13, 10}
{\it Suppose that $\A$ verifies
additional assumption (3.4)
with $p=n$ and that $\O$ is bounded and regular.
For each Radon measure $\mu$ of finite total variation,
there is a unique solution $u$ to (3.1) such that
\roster
\item"i)" the truncations
$$
T_k(u)\in W^{1,n}_0(\O)\text{ for all }k>0\,,
$$
\item"ii)" $u$ is in $\operatorname{BMO}$, and
 $$
\nabla u\in \wL{n}{\O}.
$$
\endroster
}

The regularity of $\O$ refers to the fact that the complement of $\O$
needs to be thick enough to exclude counterexamples we had in Section
2. Zhong \cite{41} formulated a weak condition for this by
requiring that the complement of $\O$ is {\it uniformly} $p$-{\it thick},
i.e\., $  \kapu(\complement\O\cap B(x,r))\approx r^{n-p}$
for all small $r>0$;
see \cite{23, 15, 32} for more information
about uniform thickness.

In fact stronger uniqueness properties than in Theorem 3.7 hold:
by using a Hodge decomposition argument
Greco, Iwaniec, and Sbordone \cite{13} proved that for $p$-Laplacian
the solution is unique
in the {\it grand Sobolev space}  $W^{1,n)}(\O)$, i.e\.
$$
u\in{\tsize\bigcap\limits_ {q<n}}
            W^{1,q}_0(\O)\quad
\text{ and }\quad\sup_{\vep>0}\vep\int_{\O}|\nabla u|^{n-\vep}\,dx<\infty\,.
$$

The regularity $\nabla u\in \wL{n}{\O}$
in Theorem 3.7  (proved in \cite{10}) is better than
$u\in W^{1,n)}(\O)$.

By using a maximal function argument similar to that introduced
by Lewis \cite{24},
Zhong \cite{41} proved a stronger uniqueness result: the solution is
unique in
$$
u\in{\tsize\bigcap\limits_ {q<n}}
            W^{1,q}_0(\O)\,.
$$
Even stronger result appeared in \cite{10}: there is a $\vep>0$ depending on
the structural assumptions and on $\O$ such that any solution in
$W^{1,n-\vep}_0(\O)$ is actually the unique solution declared in Theorem
3.7. A similar result follows from the estimates in \cite{13}.

\specialhead \noindent \boldLARGE References \endspecialhead

\ref\key 1 \by P. B\'enilan, L.  Boccardo, T. Gallou\"et, R. Gariepy, M.
 Pierre,  \& J. L. V\'azquez
\paper An $L\sp 1$-theory of existence
and uniqueness of solutions of nonlinear elliptic equations
\jour   Ann. Scuola Norm. Sup. Pisa Cl. Sci. (4)
\vol22  \yr1995\pages 241--273
\endref

\ref\key 2
\by J. Bj\"orn, P. MacManus \& N. Shanmugalingam
\paper Fat sets and pointwise boundary estimates for
$p$-harmonic functions in metric spaces
\jour J. Anal. Math.
\vol 85
\yr 2001\pages 339--369
\endref

\ref\key 3
\by L. Boccardo \& T. Gallou\"et
\paper Nonlinear elliptic equations with right-hand side measures
\jour  Comm. Partial Diff. Eq.
\vol17  \yr1992\pages 641--655
\endref

\ref\key 4
\by L. Boccardo, T. Gallou\"et, \& L. Orsina
\paper Existence and uniqueness of entropy solutions for
nonlinear elliptic equations with measure data
\jour   Ann. Inst. H. Poincar\'e Anal. Non Lin\'eaire
\vol13  \yr1996\pages  539--551
\endref

\ref\key 5
\by L. Carleson
\book Selected Problems on Exceptional Sets
\publ Van Nostrand\yr 1967
\endref

\ref
\key 6
\by G. Dal Maso \& A. Malusa
\paper Some properties of reachable solutions of
nonlinear elliptic equations with measure data
\jour   Ann. Scuola Norm. Sup. Pisa Cl. Sci. (4)
\vol 25
\yr 1997
\pages 375--396
\endref

\ref
\key 7
\by G. Dal Maso, F. Murat, L. Orsina,  \& A. Prignet
\paper Definition and existence of
renormalized solutions of elliptic equations with general measure data
\jour  C. R. Acad. Sci. Paris S\'er. I Math.
\vol 325
\yr 1997
\pages 481--486
\endref

\ref\key 8
\by G. Dal Maso, F. Murat, L. Orsina,  \& A. Prignet
\paper Renormalized solutions of elliptic equations with
general measure data
\jour Ann. Scuola Norm. Sup. Pisa Cl. Sci. (4)
\yr1999
\vol 28
\pages741--808
\endref

\ref\key 9
\by G. Dolzmann, N. Hungerb\"uhler,  \&  S. M\"uller
\paper Non-linear elliptic systems with measure-valued right hand side
\jour  Math. Z.
\vol 226
\yr  1997
\pages 545--574
\endref

\ref\key 10
\by G. Dolzmann, N. Hungerb\"uhler,  \&  S. M\"uller
\paper Uniqueness and maximal regularity for nonlinear elliptic
systems of $n$-Laplace type with measure valued right hand side
\jour J. Reine Angew. Math.  \yr 2000\vol 520
\pages1--35
\endref

\ref
\key 11
\by R. Gariepy \& W. P. Ziemer
\paper A regularity condition at the boundary for solutions of
quasilinear elliptic equations
\jour Arch. Rat. Mech. Anal.
\vol 67
\yr 1977
\pages 25--39
\endref

\ref
\key 12
\by S. Granlund, P. Lindqvist, \& O. Martio
\paper Conformally invariant variational integrals
\jour Trans. Amer. Math. Soc.\yr 1983\vol277\pages 43--73
\endref

\ref\key 13
\by L. Greco, T. Iwaniec, \& C.  Sbordone
\paper Inverting the $p$-harmonic operator
\jour  Manu\-scripta Math.
\yr 1997\vol 92\pages  249--258
\endref

\ref
\key 14
\by O. Guib\'e
\paper Remarks on the uniqueness of comparable
renormalized solutions of elliptic equations with measure data
\jour Ann. Mat. Pura Appl. (4)
\vol 180
\yr 2002
\pages 441--449
\endref

\ref\key 15
\by J. Heinonen, T. Kilpel\"ainen, \&  O. Martio
\book Nonlinear Potential Theory of Degenerate Elliptic Equations
\publ Oxford University Press\publaddr Oxford\yr 1993
\endref

\ref\key{16}
\by T. Kilpel\"ainen
\paper Nonlinear potential theory and PDEs
\jour Potential Analysis
\vol3
\yr 1994
\pages 107--118
\endref

\ref\key 17
\by  T. Kilpel\"ainen
\paper Singular solutions to $p$-Laplacian type
equations
\jour Ark. Mat.
\vol 37
 \yr1999
\pages 275--289
\endref

\ref\key 18
\by  T. Kilpel\"ainen \& J. Mal\'y
\paper Degenerate elliptic equations with measure data and nonlinear
potentials
\jour Ann. Scuola Norm. Sup. Pisa Cl. Science, Ser\. IV
\vol 19
\yr 1992
\pages591--613
\endref

\ref\key 19
\by T. Kilpel\"ainen \& J. Mal\'y
\paper The Wiener test and potential estimates for quasilinear
elliptic equations
\jour Acta Math. \vol 172\yr 1994\pages 137--161
\endref

\ref\key 20
\by T. Kilpel\"ainen \& X. Xu
\paper  On the uniqueness problem for quasilinear
elliptic equations involving measures
\jour  Rev. Mat. Iberoamer. \vol 12
\yr 1996\pages 461--475
\endref

\ref\key 21
\by T. Kilpel\"ainen \&  X. Zhong
\paper Removable sets for continuous solutions
of quasilinear elliptic equations
\jour Proc. Amer. Math. Soc.
\yr 2002
\vol 130
\pages 1681--1688
\endref

\ref
\key 22
\by D. Labutin
\paper Potential estimates for a class of fully
nonlinear elliptic equations
\jour Duke Math. J.
\vol 111
\yr 2002
\pages 1--49
\endref

\ref
\key 23
\by J. L. Lewis
\yr 1988
\paper Uniformly fat sets
\jour Trans. Amer. Math. Soc.
\vol 308
\pages 177--196
\endref

\ref
\key 24
\by J. L. Lewis
\paper On very weak solutions of certain elliptic systems
\jour Comm. Partial Diff. Eq.
\vol 18
\yr 1993
\pages 1515--1537
\endref

\ref\key 25
\by G. M. Lieberman
\paper Sharp forms of estimates for subsolutions and
supersolutions of quasilinear elliptic equations involving
measures
\jour  Comm. Partial Diff. Eq.
\vol 18
\yr 1993
\pages 1191--1212
\endref

\ref
\key 26
\by P. Lindqvist
\yr 1986
\paper On the definition and properties of
$p$-superharmonic functions
\jour J. Reine Angew. Math.
\vol 365
\pages 67--79
\endref

\ref
\key 27
\by  P. Lindqvist \& O. Martio
\paper Two theorems of N.~Wiener for solutions of quasilinear
elliptic equations
\jour Acta Math.
\vol 155
\yr 1985
\pages 153--171
\endref

\ref\key 28
\by P.L. Lions \& F. Murat
\paper Solutions renormalis\'ees d'\'equations elliptiques
non lin\'eaires
\paperinfo (to appear)
\endref

\ref\key 29
\by J. Mal\'y
\paper Examples of weak minimizers with continuous singularities
\jour  Exposition. Math.
\vol 13
\yr 1995
\pages 446--454
\endref

\ref
\key 30
\by J. Mal\'y \& W. P.  Ziemer
\book Fine Regularity of Solutions of Elliptic Partial Differential Equations
\bookinfo  Mathematical Surveys and Monographs, 51
\publ  Amer. Math. Soc.
\publaddr Providence, RI
\yr 1997
\endref

\ref
\key 31
\by V. G. Maz'ya
\paper On the continuity at a boundary point of solutions
of quasi-linear elliptic equations (English translation)
\jour Vestnik Leningrad Univ. Math.
\vol 3
\yr 1976
\pages 225--242
\finalinfo Original in {\it Vestnik Leningrad. Univ.}
{\bf 25\/} (1970), 42--55 (in Russian)
\endref

\ref\key 32
\by P. Mikkonen
\yr 1996
\paper On the Wolff potential and quasilinear elliptic equations
involving measures
\jour Ann. Acad. Sci. Fenn. Ser. A I.
Math. Dissertationes
\vol 104
\pages 1--71
\endref

\ref
\key 33
\by  A. Prignet
\paper Remarks on existence and uniqueness of solutions of elliptic
problems with right-hand side measures
\jour  Rend. Mat. Appl.
\vol 15
\yr 1995
\pages  321--337
\endref

\ref
\key 34
\by J.-M. Rakotoson
\paper Uniqueness of renormalized solutions in a $T$-set
for the $L\sp 1$-data problem and the link between
various formulations
\jour Indiana Univ. Math. J.
\vol 43
\yr 1994
\pages 685--702
\endref

\ref
\key 35
\by J.-M. Rakotoson
\paper Properties of solutions of quasilinear equations in a
$T$-set when the datum is a Radon measure
\jour Indiana Univ. Math. J.
\vol 46
\yr 1997
\pages  247--297
\endref

\ref
\key 36
\by  J.-M. Rakotoson \& W. P. Ziemer
\paper  Local behavior of solutions of quasilinear
elliptic equations with general structure
\jour  Trans. Amer. Math. Soc.
\vol 319
\yr 1990
\pages 747--764
\endref

\ref\key 37
\by J. Serrin
\yr 1964
\paper Local behavior of solutions to  quasi-linear equations
\jour Acta Math.
\vol 111\pages 247--302
\endref

\ref\key 38
\by J. Serrin
\yr 1964
\paper Pathological solutions of elliptic differential equations
\jour Ann. Scuola Norm. Sup. Pisa (3)
\vol 18\pages 385--387
\endref

\ref\key 39
\by  N.~S. Trudinger \& X.J. Wang
\paper On the Weak continuity of elliptic operators and
applications to potential theory
\jour Amer. J. Math.
\yr 2002
\vol 124
\pages 369--410
\endref

\ref
\key 40
\by  N.~S. Trudinger \& X.J. Wang
\paper Dirichlet problems for quasilinear elliptic equations with
measure data
\jour preprint
\endref

\ref\key 41
\by X. Zhong
\paper On nonhomogeneous quasilinear elliptic equations
\jour  Ann. Acad. Sci. Fenn. Math. Diss.  \vol 117 \yr1998
\endref

\enddocument